\documentclass[a4paper,10pt]{article}
\usepackage{amsfonts}
\usepackage{amscd,amsmath}
\usepackage{amssymb,amsthm}
\usepackage{euscript}
\usepackage{tabularx}
\usepackage[all]{xy}
\usepackage[T1]{fontenc}
\usepackage[T1]{fontenc}
\usepackage[utf8]{inputenc}
\usepackage[english]{babel}
\usepackage{hyperref}
\theoremstyle{plain}
\newtheorem{theorem}{Theorem}[section]

\newtheorem{lemma}[theorem]{Lemma}
\newtheorem{coro}[theorem]{Corollary}
\newtheorem{example}[theorem]{Example}

\theoremstyle{definition}

\theoremstyle{remark}

\newcommand{\lmn}[2]{\hbox to\textwidth{#1\dotfill #2}}

\newcommand{\rr}{\mathfrak{r}}
\newcommand{\Z}{\mathbb{Z}}

\newcommand{\s}{\mathfrak{s}}
\newcommand{\U}{\mathcal{U}}

\DeclareMathOperator{\lcm}{lcm}
\title{Dimension Quotients of Metabelian Lie Rings}
\author{Inder Bir S. Passi and  Thomas Sicking}

\begin{document}

\maketitle

\begin{abstract}\par\noindent
For a Lie ring $L$ over the ring of integers, we compare its lower \linebreak central series $\{\gamma_n(L)\}_{n\geq 1}$ and its dimension series  $\{\delta_n(L)\}_{n\geq 1}$ defined by setting 
$\delta_n(L)= L\cap \varpi^n(L)$, where $\varpi(L)$ is the augmentation ideal of the universal  enveloping algebra of $L$. While   $\gamma_n(L)\subseteq\delta_n(L)$ for all $n\geq 1$, the two series can differ.
In this paper it is proved that if  $L$ is a metabelian Lie ring, then   $2\delta_n(L)\subseteq\gamma_n(L)$, and $[\delta_n(L),\,L]=\gamma_{n+1}(L)$, for all $n\geq 1$.

\end{abstract}
\par\vspace{.25cm}\noindent
\par\vspace{.25cm}\noindent
\section{Introduction}
For a group $G$ and its group ring $\Z[G]$ with augmentation ideal $\varpi(G)=\linebreak \Z[G](G-1)$,
one can define the following two series of subgroups of $G$. Firstly, we have the lower central series $\{\gamma_n(G)\}_{n\geq 1}$, defined by setting
$\gamma_1(G)=G$ and $\gamma_{n+1}(G)=[G,\, \gamma_n(G)]$ for $n\geq 1$; secondly, we have the so called dimension series $\{\delta_n(G)\}_{n\geq 1}$ with   $\delta_n(G)=G\cap(1+\varpi^n(G))$.
It is easy to see that $\gamma_n(G)\subseteq\delta_n(G)$ for all $n\geq 1$, but the converse inclusion is not true in general, as first shown by  E. Rips  \cite{Rips}. Subsequently, it has been proved by J. A. Sjogren \cite{Sjogren} (see also \cite{CliffHartley}) that, for every group $G$, the quotient group $\delta_n(G)/\gamma_n(G)$ has exponent dividing $c_n$, where $$c_n=\prod_{k=1}^{n-2}b_k^{\binom{n-2}{ k}}$$ 
with  $b_k=\lcm(1,\,\ldots,\,k)$. For metabelian groups, N. Gupta \cite{Gupta:1984} (see also \cite{Gupta:1989a}) has given such  constants which are sharper than the ones obtained by Sjogren. 
\par\vspace{.25cm}
 
For Lie rings over the ring $\mathbb Z$ of integers, an analogous study has been \linebreak initiated in \cite{BartholdiPassi}. Let $L$ be a Lie ring over $\Z$ and let $\U(L)$ be its universal  \linebreak enveloping  algebra. In $L$ we  again have the 
lower central series $\{\gamma_n(L)\}_{n\geq 1}$ given by setting $\gamma_1(L)=L$ and $\gamma_{n+1}(L)=[L,\,\gamma_n(L)]$ for $n\geq 1$, and the Lie dimension  series $\{\delta_n(L)\}_{n\geq 1}$ in $L$ can be defined as follows.
Let $\varpi(L)$ be the  augmentation ideal of the universal enveloping algebra $\U(L)$ of $L$. Since  $L$ embeds into $\U(L)$ (see, e.g., \cite{Cartier}), $\varpi(L)$ is the two-sided ideal in $\U(L)$ generated by $L$.
Thus we can define the  series $\{\delta_n(L)\}_{n\geq 1}$ of Lie subrings of $L$ by setting $\delta_n(L):=L\cap\varpi^n(L)$, $n\geq 1$. It is easy to see that $\gamma_n(L)\subseteq\delta_n(L)$ for all  $n\geq 1$. It is shown in \cite{BartholdiPassi} that  $\delta_n(L) =\gamma_n(L) $ for $n\leq 3$ for every Lie ring $L$, and a counterexample to the \linebreak equality is given for the case $n=4$.  It is further  shown in \cite{BartholdiPassi} that $2\delta_4(L)\subseteq \gamma_4(L)$ \linebreak always. In view of the fact that in Lie algebras over the field of rationals the lower central and the dimension series coincide \cite{Knus}, all quotients $\delta_n(L)/\gamma_n(L)$ are  torsion abelian groups. It is shown in \cite{Sicking} that an analogue of Sjogren's result holds as well for Lie rings, i.e. $c_n\delta_n(L)\subseteq\gamma_n(L)$ for every Lie ring $L$.
\par\vspace{.25cm}

Our main result in this paper (Theorem \ref{metabelian}) is that for every metabelian Lie ring $L$, $2\delta_n(L)\subseteq \gamma_n(L)$ for all $n\geq1$. We further show that, as in the case $n=4$ given in \cite{BartholdiPassi}, there exist metabelian Lie rings $L(n)$ for all $n\geq 5$ with nontrivial $\delta_n(L(n))/\gamma_n(L(n))$.

\section{Dimension quotients of metabelian Lie rings}

If $F$ is a free Lie ring with basis  $X$, then  its  universal enveloping algebra $\U(F)$ is isomorphic to the 
free associative algebra over $\Z$ generated by $X$, and the augmentation ideal $\varpi(F)$ consists of all  polynomials whose constant term is zero.
It is shown in \cite{BartholdiPassi}, that $\gamma_n(F)=\delta_n(F)$ for all $n\geq 1$. 
\par\vspace{.25cm}
For notational convenience, for a Lie ring $L$, we denote by  $L'$  its derived subring $[L,\,L]$; thus, in particular, $L''$ denotes the second derived subring of $L$. 
\par\vspace{.25cm}
\begin{theorem}\label{metabelian}
If  $L$ is a metabelian Lie ring, then
$
 2\delta_n(L)\subseteq\gamma_n(L)
$
for all $n\geq 1$.
\end{theorem}
\par\vspace{.25cm}

To prove this theorem, it clearly suffices to consider only finitely generated Lie rings. 
\par\vspace{.25cm}
Let $L$ be a finitely generated metabelian Lie ring and let  $L\cong F/R$ be a pre-abelian  presentation of  $L$ with $F$ a free Lie ring having basis $X=\{X_1,\,\ldots\,,\,X_m\}$ and $R\supseteq F''$ an ideal  of $F$ generated by $e_1X_1+\xi_1,\,\ldots\,,\,e_mX_m+\xi_m,\,\xi_{m+1},\,\ldots$ such that $e_1,\,\ldots\,,\,e_m$ are nonnegative integers with $e_i\mid e_{i+1}$ for all $i\in\{1,\ldots,m-1\}$, and $\xi_i\in F'$, the derived subring of $F$. It is easy to see that every finitely generated Lie ring has such a pre-abelian presentation.  Denote by $\rr$ the ideal of $\U (F)$ generated by $R$. We then have $\U(L)\cong \U (F)/\rr$.
Furthermore, we have $\gamma_n(L)\cong (\gamma_n(F)+R)/R$ and $\delta_n(L)\cong (F\cap(\varpi(F)^n+\rr))/R$ for all $n\geq 1$. Observe that 
\begin{equation}\label{R=FR}
\rr=\varpi(F)\rr+R;
\end{equation}
thus, for any element $f\in F\,\cap\,(\varpi^n+\rr)$, there exists an element $f^\prime\in\linebreak  F\cap(\varpi^n+\varpi\rr)$, such that
$f^\prime\equiv f\mod R$.\par\vspace{.25cm}
Set $S=F'+R$, and denote by $\s$ the two-sided ideal of $\U(F)$ generated by $S$. 
We begin by identifying the Lie ideals $$M:=F\cap\varpi(F)\s\quad \text{and}\quad F\cap(\varpi(F)^n+\varpi(F)\s).$$
\par\vspace{.25cm}

\begin{lemma}\label{fs}
With the notation as above, we have 
 \newcounter{met1} 
 \begin{list}{$(\roman{met1})$}{\usecounter{met1}}
   \item $M=\langle e_i[X_i,\,X_j]\,|\, m\,\geq\,i\,>\,j\geq 1\rangle+[F^{\prime},\,S]$, 
   \item $[M,\, F]\subseteq [F^\prime,\, S]$,
   \item $F\cap(\varpi(F)^n+\varpi(F)\s)=\gamma_n(F)+M$ for all $n\geq 1$.
 \end{list}
\end{lemma}
\par\vspace{.25cm}

 \begin{proof} $(i)$ 
  It is clear that $[F^\prime,\,S]$ is contained in $F\cap\varpi(F)\s$. Also, for $i>j$ we have $e_j\mid e_i$, and thus $e_i[X_i,\,X_j]=X_ie_iX_j-X_je_iX_i\in F\cap\varpi(F)\s$.
  \par\vspace{.25cm}
  To check the reverse inclusion,  note first that  $F$ acts on itself via the adjoint action:
  $$x^y:=ad_y(x)=[x,\,y] \ \text{for}\ x,\,y\in F.$$  
  This action extends to an action of the universal enveloping algebra  $\U(F)$ on $F$, and for an associative monomial $w=X_{i_1}\cdots X_{i_d}\in \U(F)$
  and $x\in F$, we have $$x^w=[x,\,X_{i_1},\,\ldots\,,\,X_{i_d}],$$ where the Lie  commutator is left normed. 
\par\vspace{.25cm}
Now  suppose  $v\in F\cap\varpi(F)\s$. Then $v\in  F\cap \varpi(F)^2=F^\prime$.  Therefore, 
  \[
   v\equiv\sum_{1\leq j\leq m-1}\sum_{m\,\geq\,i\,>\,j}[X_i,X_j]^{u_{ij}} \mod F^{\prime\prime},
  \]
with  $u_{ij}\in \mathcal U(F)$ involving only the variables $X_j,\,\ldots\,,\,X_m$.  Let  $$v_j:=\sum_{m\,\geq\,i\,>\,j}[X_i,\,X_j]^{u_{ij}}.$$ For $j=1,\,\ldots\,,\, m-1$, define
endomorphisms $\theta_j$ of $\U(F)$ by setting $\theta_j(X_i)=0$ for $i<j$ and $\theta_j(X_i)=X_i$ for $i\geq j$. Since  both $\varpi(F)$
and $\s$ are  invariant under each  $\theta_j$, successive application of the
endomorphisms $\theta_{m-1},\ldots,\theta_1$ shows that $v_j\in F\cap\varpi(F)\s$ for
$1\leq j\leq m-1$.  
Write 
$$u_{ij}=n_{ij}+u^\prime_{ij}\mbox{ with } n_{ij}\in\Z,\,u^\prime_{ij}\in\varpi(F)$$
so that 
\begin{equation}\label{m1}
 v_j=\sum_{m\,\geq\,i\,>\,j}n_{ij}[X_i,X_j]+\sum_{m\,\geq\,i\,>\,j}[X_i,X_j]^{u^\prime_{ij}}.
\end{equation}
Note that, modulo $\varpi(F)\s$, the action of $\U(F)$ is just right multiplication. Thus we have, modulo $\varpi(F)\s$,
\begin{equation}\label{vjzero}
 0\equiv v_j\equiv\sum_{m\,\geq\,i\,>\,j}n_{ij}(X_iX_j-X_jX_i)+\sum_{m\,\geq\,i\,>\,j}(X_iX_j-X_jX_i)u^\prime_{ij}.
\end{equation}
The right coefficient of $X_i$, $(m\,\geq\,i\,>\,j)$, in this expression is $X_ju_{ij}$ and must be in $\s$. As $u_{ij}$ only involves the variables $X_j,\,\ldots\,,\,X_m$,
it follows that $u_{ij}$ is divisible by $e_j$, so that $n_{ij}$ is divisible by $e_j$ and $$u_{ij}^\prime=e_jf_{ij}.$$
Now the right coefficient of $X_j$ is $\sum_{m\,\geq\,i\,>\,j}(-n_{ij}X_i-X_i u^\prime_{ij})$ and it must also lie in $\s$. Thus $e_i\mid n_{ij}$ for each $i>j$,
as the only terms of degree one in $\s$ are linear combinations of  $e_kX_k$, $1\leq k\leq m$. Thus it follows that  $$\sum_{m\,\geq\,i\,>\,j}X_i u_{ij}^\prime \in\s.$$

\par\vspace{.25cm}Denote by $\mathfrak{a}$ the ideal of $\U (F)$ generated by $F^\prime$. As $\mathfrak{a}$ acts trivially on $[X_i,X_j]$ modulo $F^{\prime\prime}$,
we can change $u_{ij}^\prime$ modulo $\mathfrak{a}$, and thus assume, without loss of generality,  that $u_{ij}^\prime$ and hence $f_{ij}$ is a linear combination of ordered monomials, i.e., monomials of the
form $X_{k_1}\cdots X_{k_t}$ with $j\leq k_1\leq \ldots\leq k_t\leq m$. Thus we can write
\[
 f_{ij}=\sum_{k=j}^m X_k g_{ijk},
\]
with  $g_{ijk}$ involving only the variables $X_k,\,\ldots\,,\, X_m$. Thus  we have 
\[
 \sum_{m\,\geq\,i\,>\,j}X_i u_{ij}^\prime=\sum_{m\,\geq\,i\,>\,j}\sum_{k=j}^{m}e_j X_iX_k g_{ijk}\in \s.
\]
As $e_jX_iX_jg_{ijj}\in\s$, we also have $\sum_{m\,\geq\,i\,>\,j}\sum_{k=j+1}^m e_j X_iX_k g_{ijk}\in \s$. \linebreak Successively applying $\theta_{j+1},\,\ldots,\,\theta_{i-1}$
it follows inductively that $e_k$ divides\linebreak 
$\sum_{m\,\geq\,i\,>\,j}e_jX_iX_kg_{ijk}$ for each $k$ with $j\leq k< i$. Thus we can write
\[
 \sum_{m\,\geq\,i\,>\,j} X_i u_{ij}^\prime=\sum_{m\,\geq\,i\,>\,j}\left(\sum_{k=j}^{i-1}e_k X_iX_k g_{ijk}^\prime+e_i h_{ij}\right),
\]
where $h_{ij}\in\varpi(F)$ involves only the variables $X_i,\,\ldots,\,X_m$.
Consequently,
\[
 \sum_{m\,\geq\,i\,>\,j}[X_i,X_j]^{u_{ij}^\prime}=\sum_{m\,\geq\,i\,>\,j}\sum_{k=j}^{i-1} e_k[X_i,\,X_j,\,X_k]^{g_{ijk}^\prime}+\sum_{m\,\geq\,i\,>\,j}e_i[X_i,X_j]^{h_{ij}}.
\]
Now $e_k[X_i,\,X_j,\,X_k]=[X_i,\,X_j,\,e_kX_k]\in[F^\prime,\,S]$.
For $k\geq i$, using  Jacobi \linebreak identity, we have $e_i[X_i,\,X_j,\,X_k]\in [F^\prime,\,S]$. Since $h_{ij}\in \varpi(F) $ is a polynomial only in the variables $X_k$ with  $k\geq i$, it follows that  $e_i[X_i,\,X_j]^{h_{ij}}\in [F^\prime,\,S]$. Since $F^{\prime\prime}$ is contained in $[F^\prime,\,S]$, the proof of (i) is complete.
\par\vspace{.5cm}\noindent 
$(ii)$ follows immediately from (i) in view of the Jacobi identity.

\par\vspace{.5cm}\noindent 
$(iii)$ Clearly, $\gamma_n(F) +M\subseteq F\cap(\varpi^n(F)+\varpi(F)\mathfrak s)$ for all $n\geq 1$. For $n=1,\,2$, the reverse inclusion  holds trivially. 
Suppose  $n\geq 3$ and    $v\in F\cap(\varpi(F)^n+\varpi(F)\s)$.  Then, proceeding as in case (i) above,  we have 
\[
   v\equiv\sum_{m\,\geq\,i\,>\,j\geq\, 1}[X_i,\,X_j]^{u_{ij}}=\sum_{j=1}^{m-1} v_j \mod ( F^{\prime\prime}+\gamma_n(F)).
\]
and, modulo $\varpi^n(F)+\varpi(F)\s$, 
\begin{equation}\label{vjzero1}
 0\equiv v_j\equiv\sum_{m\,\geq\,i\,>\,j}n_{ij}(X_iX_j-X_jX_i)+\sum_{m\,\geq\,i\,>\,j}(X_iX_j-X_jX_i)u^\prime_{ij}.
\end{equation}
Therefore,  the right coefficient of $X_j$  must be in $\varpi^{n-1}(F)+\s$.  Hence $n_{ij}$ must be divisible by $e_i$ and, modulo $\varpi^{n-2}(F)$,  
$u_{ij}^\prime$ must be divisible by $e_j$. Writing $u_{ij}^\prime=e_jf_{ij}+p_{ij}$ with $p_{ij}\in\varpi^{n-2}(F)$ and decomposing $f_{ij}$ as in the above proof of (i), we have
\begin{multline*}
 \sum_{m\,\geq\,i\,>\,j}[X_i,X_j]^{u_{ij}^\prime}=\sum_{m\,\geq\,i\,>\,j}\sum_{k=j}^{i-1} e_k[X_i,X_j,X_k]^{g_{ijk}^\prime}\\+\sum_{m\,\geq\,i\,>\,j}e_i[X_i,X_j]^{h_{ij}}+\sum_{m\geq i>j}[X_i,\,X_j]^{p_{ij}}\in M+\gamma_n(F).
\end{multline*}
Hence $v\in \gamma_n(F)+M$.
 \end{proof}

\par\vspace{.25cm}

We immediately have the following corollary.

\par\vspace{.25cm}
\begin{coro}
 If $L$ is a metabelian Lie ring, then $
  [\delta_n(L),\,L]=\gamma_{n+1}(L)
 $
 for all $n\geq 1$.
 \begin{proof}
 Since $\gamma_n(L)\subseteq \delta_n(L)$, we clearly have that $\gamma_{n+1}(L)\subseteq [\delta_n(L),\,L]$ for all $n\geq 1$.
\par\vspace{.25cm}\noindent 
 For the reverse inclusion, let $L=F/R$ be a free presentation of a metabelian Lie ring $L$ and let $S=R+F^\prime$. Then $[F^{\prime},\,S]\subseteq R$, and, by Lemma \ref{fs}, we have 
 \[
  [F\cap(\varpi^n(F)+\varpi(F)\rr),\,F]\subseteq [F\cap(\varpi^n(F)+\varpi(F)\s),\,F]\subseteq \gamma_{n+1}(F)+R.
 \]
 
Hence it follows that $[F\cap(\varpi^n(F)+\rr),\,F]\subseteq \gamma_{n+1}(F)+R$, and the Corrolary is proved.
 
 \end{proof}

\end{coro}

\subsection{Proof of the Theorem \ref{metabelian}}

Let $L$ be a finitely generated Lie ring and $L=F/R$ a pre-abelian presentation of $L$ as set out at the beginning of Section 2.  Suppose   $v\in F\cap(\varpi^n(F)+\rr)$. Then there exists  $r\in R$ such that $w:=v+r\in F\cap(\varpi^n(F)+\varpi(F)\rr)\subseteq F\cap(\varpi(F)^n+\varpi(F)\s)$. By  Lemma \ref{fs}, we then have 
\begin{equation}\label{met}
 w\equiv \sum_j \sum_{m\,\geq\,i\,>\,j} d_{ij}[X_i,\,X_j]+\sum_k e_k[X_k,\,\eta_k] \mod (\gamma_n(F)+F^{\prime\prime})
\end{equation}
with $d_{ij}\in\Z$ such that $e_i\mid d_{ij}$, and $\eta_k\in F^\prime$. In $\U(F)$ this congruence yields
\begin{align*}
 w\equiv \sum_j\sum_{m\,\geq\,i\,>\,j}d_{ij}(X_iX_j-X_jX_i)+\sum_k e_k X_k \eta_k \mod \varpi^2(F)\s + \varpi^n(F). 
\end{align*}
Since $w\in \varpi(F)^n+\varpi(F)\rr$, it follows that 
$$ \sum_k X_k(y_k+e_k \eta_k)\equiv 0 \mod (\varpi^2(F)\s + \varpi^n(F) + \varpi(F)\rr),$$
where \[
 y_k=-\sum_{m\geq i>k}d_{ik}X_i+\sum_{1\leq j<k}d_{kj}X_j\in F.
\]
Thus,  for every $k$,  we have $y_k+e_k\eta_k\in \varpi(F)\s+\varpi^{n-1}(F)+\rr$, and  there exists $r_k\in\rr\cap F=R$ with
$y_k+e_k\eta_k+r_k\in \varpi(F)\s+\varpi^{n-1}(F)$. Since  $[F^\prime,\,S]\subseteq S^\prime\subseteq R$, by Lemma \ref{fs},  we deduce that
\[
 [X_k,\,y_k+e_k\eta_k+r_k]\in \gamma_n(F)+[F^\prime,\,S]\subseteq \gamma_n(F)+R.
\]
As $[X_k,\,r_k]\in R$ and $[X_k,\,e_k\eta_k]=[e_k X_k,\,\eta_k]\in R$,
we conclude that $$[X_k,\,y_k]\in R+\gamma_n(F)\ \text{ for every}\  k.$$ Finally, the congruence (\ref{met}) reduces, modulo $\gamma_n(F)+R$, to yield 
\begin{equation}\label{met1}
 w\equiv \sum_j\sum_{m\,\geq\,i\,>\,j}d_{ij}[X_i,\,X_j]=\sum_k X_ky_k.
\end{equation}
Recall that the map given by $u\mapsto -u$ on $F$ induces  an anti-automorphism on $\U(F)$. Applying this anti-automorphism  to the above equation yields
\begin{equation}\label{met2}
 -w\equiv \sum_k y_kX_k.
\end{equation}
From the congruences (\ref{met1}) and (\ref{met2}), we have 
\[
 2w=\sum_k X_ky_k-y_kX_k=\sum_k[X_k,\,y_k]\in \gamma_n(F)+R,
\]
and the proof is complete.  $\Box$
\par\vspace{.5cm}
Finally, following Gupta's construction of counter examples to the dimension conjecture for groups, as given in \cite{Gupta:1989}, and the construction of a Lie ring $L$ with $\delta_4(L)\neq \gamma_4(L)$ in \cite{BartholdiPassi}, we observe that analogously, for each $n\geq 5$, there exists a metabelian Lie ring $L(n)$ with nontrivial $\delta_{2n-4}(L(n))$ and $\gamma_n(L(n))=0$.
\par\vspace{.25cm}
\begin{example}
Let $X=\{r,a,b,c\}$ and let $F$ be the free Lie ring generated by $X$. Set $$x_0=y_0=z_0=r,$$ $$x_i=[x_{i-1},\,a],\, y_i=[y_{i-1},\,b],\, z_i=[z_{i-1},\,c] \mbox{ for } i\geq 1.$$ Let $R(n),\ n\geq 4,$ be the Lie ideal of $F$ generated by the following elements:
\par\vspace{.25cm}\noindent
$2^{2n-1}r,$\\
 $2^{n+2}a-4y_{n-3}-2z_{n-3},\, 2^nb+4x_{n-3}-z_{n-3},\, 2^{n-2}c+2x_{n-3}+y_{n-3},$\\
$z_{n-2}-4y_{n-2},\, y_{n-2}-4x_{n-2},$\\
$x_{n-1},\,y_{n-1},\,z_{n-1},$\\
$[a,\,b,\,u],\,[a,\,c,\,u],\,[b,\,c,\,u]$ for all $u\in F,$\\
$[x_i,\,b],\,[x_i,\,c],\,[y_i,\,a],\,[y_i,\,c],\,[z_i,\,a],\,[z_i,\,b]$, $i\geq 1,$\\
$[x_i,\,x_j],\,[x_i,\,y_j],\,[x_i,\,z_j],\,[y_i,\,y_j],\,[y_i,\,z_j],\,[z_i,\,z_j]$, $i,\,j\geq 0.$
\par\vspace{.25cm}\noindent
Then $L(n):=F/R(n)$ is easily seen to be metabelian with $\gamma_n(L(n))=0$. Let
\[
 g:=2^{2n-1}[a,\,b]+2^{2n-2}[a,\,c]+2^{2n-3}[b,\,c].
\]
\end{example}
\par\vspace{.5cm}\noindent
It is straightforward to check that $g\in\varpi^{2n-4}(L)$, i.e., $g\in\delta_{2n-4}(L(n))$, and that $$g=2^{n+1}x_{n-2}=2^{n-3}z_{n-2}=2^{2n-1}[a,\,b].$$ That $g\neq 0$ is easily seen by comparing (in $F$) the element  $2^{n-3}z_{n-2}+\gamma_n(F)$ with the expression for an arbitrary element of $(R(n)+\gamma_n(F))/\gamma_n(F)$ in terms of the Witt basis \cite{Witt} for $F/\gamma_n(F)$ as a free abelian group. 

\par\vspace{.25cm}\noindent

\section*{Acknowledgement}
The second author would like to thank Indian Institute of Science Education and Research Mohali for the warm hospitality during his stay at the Institute. He is also thankful to Mathematisches Institut der Universit\"{a}t G\"{o}ttingen for financial support.

\par\vspace{1cm}
\noindent
Inder Bir S. Passi\\
Centre for Advanced Study in Mathematics\\
Panjab University, Chandigarh, India\\
and\\
Indian Institute of Science Education and Research\\
Mohali, India\\
email: ibspassi@yahoo.co.in

\par\vspace{.5cm}\noindent
Thomas Sicking\\
Mathematisches Institut\\
Georg-August-Universit\"at\\ G\"ottingen, Germany\\
email: thomas.sicking@yahoo.de

\begin{thebibliography}{99}

\bibitem{BartholdiPassi} Laurent Bartholdi and Inder Bir S. Passi: Lie dimension subrings, {\it Int. J. Algebra Comput.} {\bf 25},
1301--1325
\bibitem{Cartier}  Pierre Cartier: Remarques sur le th\'{e}or\`{e}me de Birkhoff-Witt, {\it Ann. Scuola Norm. Sup. Pisa
(3)} {\bf 12} (1958), 1--4 (French).
\bibitem{Gupta:1984} Narain Gupta: Sjogren's theorem for dimension subgrousps - the metabelian case. In: Combinatorial Group Theory and Topology, {\it Annals of Math. Study.} {\bf 111}, 197--211 (1987).
\bibitem{Gupta:1989a} Gupta, Narain: A solution of the dimension subgroup problem. {\it J. Algebra} {\bf 138} (1991), no. 2, 479–490.

\bibitem{Gupta:1989} Narain Gupta: The dimension subgroup conjecture, {\it Bull. London Math. Soc.}, {\bf 22}, 453--456 (1990).
\bibitem{CliffHartley} G. Cliff and B. Hartley: Sjogren's theorem on dimension subgroups,\linebreak {\it J. Pure Appl. Algebra} {\bf 47} (1987),
231--242.
\bibitem{Knus}Max A. Knus: On the enveloping   algebra and the descending central series of a Lie algebra, {\it J. Algebra} {\bf 12} (1969), 335 - 338. 
\bibitem{Rips} E. Rips: On the fourth integer dimension subgroup, {\it Israel J. Math.} {\bf 12} (1972), 342--346.
\bibitem{Sicking} Thomas Sicking, Lie Dimension Quotients, {\it unpublished}.
\bibitem{Sjogren} J. A. Sjogren: Dimension and lower central subgroups, {\it J. Pure Appl. \linebreak Algebra} {\bf 14} (1979), 175--194.
\bibitem{Witt} Ernst Witt: Treue Darstellung Liescher Ringe. {\it J. Reine Angew. Math.} {\bf 177} (1937), 152–160. 



\end{thebibliography}
\end{document}